\documentclass[preprint,11pt]{imsart}
\RequirePackage[T1]{fontenc}
\usepackage{color}
\usepackage[latin1]{inputenc}
\usepackage[T1]{fontenc}
\RequirePackage{amsfonts,amsthm,amsmath}
\RequirePackage[colorlinks,citecolor=blue,urlcolor=blue]{hyperref}
\usepackage{a4wide}
\usepackage{here}
\usepackage{float}
\startlocaldefs
\theoremstyle{plain}

\endlocaldefs
\newtheorem{theo}{Theorem}

\newtheorem{prop}[theo]{Proposition}

\newcommand{\zak}{\nobreak \ifvmode \relax \else
     \ifdim\lastskip<1.5em \hskip-\lastskip
     \hskip1.5em plus0em minus0.5em \fi \nobreak
     $\Box$\fi\\}

\hypersetup{pdfstartview=FitH}
\def\eps{{\epsilon}}

\def\e{{\mathbb E}}
\def\1{{\mathbb 1}}

\def\cov{\mathop{\rm Cov}\nolimits}    
\usepackage{graphicx}
\begin{document}
\doublespacing
\begin{frontmatter}
\title{Distribution of residual autocorrelations for multiplicative seasonal ARMA models with uncorrelated but non-independent error terms\\[0.5cm]
}
\runtitle{Diagnostic checking in weak SARMA models}
%
\begin{aug}
\author{\fnms{\Large{Yacouba}}
\snm{\Large{Boubacar Ma\"{\i}nassara}}
\corref{}\thanksref{t2}
\ead[label=e2]{yacouba.boubacar$\_$mainassara@univ-fcomte.fr}} \and \
\author{\fnms{\Large{Abdoulkarim}}
\snm{\Large{Ilmi Amir}}
\ead[label=e1]{osmanilmi@hotmail.com}}
\thankstext{t2}{Corresponding author}
\runauthor{Y. Boubacar Ma\"{\i}nassara and A. Ilmi Amir}
\affiliation{Universit\'e Bourgogne Franche-Comt\'e}
\address{\hspace*{0cm}\\
Universit\'e Bourgogne Franche-Comt\'e, \\
Laboratoire de math\'{e}matiques de Besan\c{c}on, \\ UMR CNRS 6623, \\
16 route de Gray, \\ 25030 Besan\c{c}on, France.\\[0.2cm]
\printead{e2,e1}}
\end{aug}
\vspace{0.5cm}
\begin{abstract}
In this paper we consider portmanteau tests for testing the adequacy of multiplicative seasonal autoregressive
moving-average (SARMA) models under the assumption that the errors are uncorrelated but not necessarily independent.
We relax the standard independence assumption on the error term in order to extend the range of application of the SARMA models.
We study the asymptotic distributions of residual and normalized residual empirical autocovariances and autocorrelations under
weak assumptions on the noise.  We establish the asymptotic behaviour of the proposed statistics.
 A set of Monte Carlo experiments and an application to monthly mean total sunspot number are presented.
\end{abstract}
\begin{keyword}[class=AMS]
\kwd[Primary ]{62M10}
\kwd{62F03}
\kwd{62F05}
\kwd[; secondary ]{91B84}
\kwd{62P05}
\end{keyword}
\begin{keyword}
Goodness-of-fit test, quasi-maximum likelihood estimation, Box-Pierce and Ljung-Box portmanteau tests,
residual autocorrelation, self-normalization, weak SARMA models
\end{keyword}
\end{frontmatter}
%
\section{Introduction}
The multiplicative seasonal autoregressive  moving average (SARMA) model of order $(p,q)(P,Q)_s$ for the univariate time series $X=(X_t)_{t\in\mathbb Z}$, is defined by
\begin{equation}
\label{ARMA}
a_{\theta_0}(L) \mathbf{a}_{\theta_0}(L){X}_t=b_{\theta_0}(L) \mathbf{b}_{\theta_0}(L){\epsilon}_t,\quad
\forall t\in \mathbb Z,
\end{equation}
where $\theta_0=(a_{01},\dots{},a_{0p},b_{01},\dots{},b_{0q},\mathbf{a}_{01},\dots{},\mathbf{a}_{0P},\mathbf{b}_{01},\dots{},\mathbf{b}_{0Q})'$ and
where the nonseasonal AR and MA operators are defined by
$a_{\theta_0}(L)=1-\sum_{i=1}^pa_{0i}L^i$ and $b_{\theta_0}(L)=1-\sum_{i=1}^qb_{0i}L^i$,
respectively, while the seasonal AR and MA operators are given by $\mathbf{a}_{\theta_0}(L)=1-\sum_{i=1}^P\mathbf{a}_{0i}L^{si}$ and $\mathbf{b}_{\theta_0}(L)=1-\sum_{i=1}^Q\mathbf{b}_{0i}L^{si}$,  respectively,  $s$ denotes the length of the seasonal period  and $L$ stands for the backshift operator. It is assumed that the model defined by \eqref{ARMA} is stationary, invertible and not redundant.
Without loss of generality, we also assume that $a_{0p}^2+b_{0q}^2+\mathbf{a}_{0P}^2+\mathbf{b}_{0Q}^2\neq 0$ (by convention ${a}_{00}={b}_{00}=\mathbf{a}_{00} =\mathbf{b}_{00}=1$).

In the standard  situation ${\epsilon}=({\epsilon}_t)_{t\in\mathbb Z}$ is assumed to be a sequence of independent and identically distributed (iid for short)  random variables with zero mean and common variance. In this standard framework, $(\epsilon_t)$ is said to be a \emph{strong
white noise} and the representation (\ref{ARMA})  is called a strong   SARMA$(p,q)(P,Q)_s$ process. In contrast with this previous definition, the representation  (\ref{ARMA}) is said to be a weak SARMA$(p,q)(P,Q)_s$ if the noise process $(\epsilon_t)$  is a \emph{weak white noise}, that is, if it satisfies
\begin{itemize}
\item[\hspace*{1em} {\bf (A0):}]
\hspace*{1em} $\mathbb{E}({{\epsilon}}_t)=0$,
$\mbox{Var}\left({{\epsilon}}_t\right)=\sigma^2_0$ and
$\mbox{Cov}\left({{\epsilon}}_t,{{\epsilon}}_{t-h}\right)=0$ for all
$t\in\mathbb{Z}$ and all $h\neq 0$.
\end{itemize}
 A strong white noise is obviously  a weak white
noise, because independence entails uncorrelatedness, but the
reverse is not true. It is clear from these definitions that the following inclusions hold:
$$\left\{\text{strong SARMA}(p,q)(P,Q)_s \right\}\subset\left\{\text{weak SARMA}(p,q)(P,Q)_s\right\}.$$
After estimating the SARMA process, the next important step in the modeling consists in checking if the estimated model fits satisfactorily the data. Thus, under 
the null hypothesis that the model has been correctly identified, the residuals ($\hat{\epsilon}_t$)
are approximately a white noise. This adequacy checking step validates or invalidates the choice of the orders $(p,q)$ and $(P,Q)_s$.

Based on the residual empirical autocorrelations $\hat\rho(h)=\sum_{t=1+h}^n\hat{\epsilon}_t\hat{\epsilon}_{t-h}/\sum_{t=1}^n\hat{\epsilon}_t^2$, where $n$ is the length of the series, \cite{bp70} have 
proposed  a goodness-of-fit test, the so-called portmanteau test, for strong ARMA models. A modification of their test has been proposed by \cite{lb}.  It is nowadays one of the most popular diagnostic checking tools in ARMA modeling of time series. These tests are defined by
\begin{equation}\label{bp}
 Q^{\textsc{bp}}_m=n\sum_{h=1}^m\hat\rho^2(h) \text{ and }   {Q}_m^{\textsc{lb}}=n(n+2)\sum_{h=1}^m\frac{\hat\rho^2(h)}{n-h},
\end{equation}
where $m$ is a fixed integer.
The  statistic ${Q}_m^{\textsc{lb}}$ has the same asymptotic chi-squared distribution as $Q_m^{\mathrm{BP}}$ and has the reputation of doing better for small or medium sized sample (see \cite{lb}).
For weak ARMA models, \cite{frz} show that the asymptotic distributions of the statistics defined
in \eqref{bp} are no longer chi-square distributions but a mixture of chi-squared distributions,
weighted by eigenvalues of the asymptotic covariance matrix of the vector of autocorrelations. Recently, \cite{BMS2018} proposed an alternative method based on a self-normalization approach to construct a new test statistic which is asymptotically distribution-free under the null hypothesis.

In many situations, these tests are implemented to check the lack of fit of SARMA models.
However, the traditional methodology of Box and Jenkins cannot be extended to the case of SARMA
models when $s>1$ because of the multiplicative  contraints on the parameters. This standard methodology needs to be adapted to take into account the possible lack of independence of the errors terms. See, for instance, \cite{Duchesne2007} and \cite{McLeod1978} who considered serial correlation testing in multiplicative seasonal univariate time series models. Duchesne (see \cite{Duchesne2007}) proposed his test statistic based on a kernel-based spectral density estimator, whose weighting scheme is more adapted to autocorrelations associated to seasonal lags.
The standard tests procedure consist in rejecting the null hypothesis of a SARMA$(p,q)(P,Q)_s$ model if the statistics (\ref{bp})  are larger than a certain quantile of a chi-squared distribution with $m-(p+q+P+Q)>0$ degrees of freedom. Consequently, these standard tests are not applicable for $m\leq p+q+P+Q$.

The works on the  portmanteau statistic of SARMA models are generally performed under the assumption that the errors $\epsilon_t$ are independent.
This independence assumption is often considered too restrictive by
practitioners. It precludes conditional heteroscedasticity and/or other forms of nonlinearity (see \cite{fz05}, for a review on weak univariate ARMA models).
In this framework, we relax the standard independence assumption on the error term in order to be able to cover SARMA representations of general nonlinear models.
For the asymptotic theory of weak SARMA, notable exception are \cite{BMIA2018} where  the consistency and the asymptotic
normality of the quasi-maximum likelihood estimator (QMLE) for weak multivariate SARMA models are studied. They also study a particular case of the asymptotic distributions of
residual autocovariances and autocorrelations at the seasonal lags $1s, 2s, 3s,\dots, ms$.

This paper is devoted to the problem of the validation step of weak SARMA representations.
We consider  portmanteau test statistics  based on the residual empirical autocorrelations
but not necessarily  at multiple lags $s$ as in \cite{BMIA2018}.
For such models, we show that the asymptotic distributions of the statistics defined in (\ref{bp}) are no
longer chi-square distributions but a mixture of chi-squared distributions, weighted by eigenvalues of  the
asymptotic covariance matrix of the vector of autocorrelations as in \cite{BMIA2018}. We also proposed another modified statistics based on a self-normalization approach  which are asymptotically distribution-free under the null hypothesis and generalize the result of \cite{BMS2018}.

In Monte Carlo experiments, we illustrate that the proposed test statistics  have reasonable finite sample performance. Under nonindependent errors,
it appears that the standard test statistics are generally non reliable, overrejecting  severely, while the proposed tests statistics offer satisfactory levels.
Even for independent errors, they seem  preferable to the standard ones, when the number $m$ of autocorrelations is small. Moreover, the error of first kind is well controlled.
Contrarily to the standard tests \eqref{bp}, the proposed tests can be used safely for $m$ small (see for instance Figure~\ref{fig}). For all these above reasons, we think that the modified
versions that we propose in this paper are preferable to the standard ones for diagnosing SARMA models under nonindependent errors.
Other contribution is  to improve the results concerning the statistical analysis of weak SARMA models by considering the adequacy problem.

The article is organized as follows. In the next section, we briefly recall the results on the QMLE asymptotic distribution obtained by \cite{BMIA2018} when $(\epsilon_t)$ satisfies mild mixing
assumptions.
We study the asymptotic behaviour of the residuals autocovariances and autocorrelations under weak assumptions on the noise in Section \ref{diagnostic1}. It is also shown
how the standard portmanteau tests \eqref{bp} must be adapted in the case of multiplicative seasonal ARMA models with non-independent innovations.
In Section \ref{diagnostic2} we derive the asymptotic distribution of residuals autocovariances and autocorrelations using self-normalization approach
and we establish the asymptotic behaviour of the proposed statistics. Section~\ref{ne} proposes numerical illustrations and an illustrative application on real data.  We provide a conclusion in Section \ref{conclusion}. The technical proofs are relegated to
the appendix.
\section{Estimating weak SARMA models}
\label{model}
In this section, we recall the results on the QMLE asymptotic distribution obtained by \cite{BMIA2018} when $(\epsilon_t)$ satisfies mild mixing
assumptions in  order to have a self-containing paper.

The unknown parameter of interest $\theta_0$ is supposed to belong to the parameter space
\begin{align*}
\Theta & =  \Big
\{\theta=(a_{1},\dots{},a_{p},b_{1},\dots{},b_{q},\mathbf{a}_{1},\dots{},\mathbf{a}_{P},\mathbf{b}_{1},\dots{},\mathbf{b}_{Q})'
\in{\Bbb R}^{k_0},\text{ where } k_0=p+q+P+Q,\\
& a_{\theta}(z)=1-\sum_{i=1}^{p}a_{i}z^{i},\, \mathbf{a}_{\theta}(z)=1-\sum_{j=1}^{P}\mathbf{a}_{j}z^{sj},\, b_{\theta}(z)=1-\sum_{i=1}^{p}b_{i}z^{i}\text{
and } \mathbf{b}_{\theta}(z)=1-\sum_{j=1}^{Q}\mathbf{b}_{j}z^{sj} \\
& \hspace{0.5cm} \text{ have all
their zeros outside the unit disk and have no zero in common}\Big \}\ .
\end{align*}
To ensure the asymptotic theory of the QMLE,
we assume that the parametrization satisfies the following smoothness conditions.
Without loss of generality, we may assume that $\Theta$ is compact.
\begin{itemize}
\item[\hspace*{1em} {\bf (A1):}]
\hspace*{1em}  {The process $\epsilon=(\epsilon_t)_{t\in\mathbb Z}$ is  ergodic and strictly  stationary}.
\end{itemize}
For the asymptotic normality of the QMLE, additional assumptions are required. It is necessary to assume
that $\theta_0$ is not on the boundary of the parameter space
${\Theta}$.
\begin{itemize}
\item[\hspace*{1em} {\bf (A2):}]
\hspace*{1em} We have
$\theta_0\in\stackrel{\circ}{{\Theta}}$, where $\stackrel{\circ}{{\Theta}}$ denotes the interior of ${\Theta}$.
\end{itemize}
To control the serial dependence of the stationary process $(\epsilon_t)$, we introduce the strong mixing coefficients $\alpha_{\epsilon}(h)$ defined by
$$\alpha_{\epsilon}\left(h\right)=\sup_{A\in\mathcal F^t_{-\infty},B\in\mathcal F_{t+h}^{+\infty}}\left|\mathbb{P}\left(A\cap
B\right)-\mathbb{P}(A)\mathbb{P}(B)\right| ,$$
where $\mathcal F_{-\infty}^t=\sigma (\epsilon_u, u\leq t )$ and $\mathcal F_{t+h}^{+\infty}=\sigma (\epsilon_u, u\geq t+h )$.
We use $|\cdot|$ to denote the Euclidian norm $|z|=\left(\sum_{i=1}^\tau z_i^2\right)^{1/2}$ of a column vector $z=(z_1,\dots,z_\tau)'$.
We will make an integrability assumption on the moment of the noise and a summability condition on the strong mixing coefficients $(\alpha_{\eps}(k))_{k\ge 0}$.
\begin{itemize}
\item[\hspace*{1em} {\bf (A3):}]
\hspace*{1em}
$\text{We have }\mathbb{E}\big|{\epsilon}_t|^{4+2\nu}<\infty\text{ and }
\sum_{k=0}^{\infty}\left\{\alpha_{{\epsilon}}(k)\right\}^{\frac{\nu}{2+\nu}}<\infty\text{ for some } \nu>0.$
\end{itemize}
Assumption {\bf (A3)} from \cite{fz98,fz05} is a technical condition for proving the asymptotic theory of the QMLE.
The integrability assumption on the moment of the noise is not very restrictive in this framework because
the innovation process $(\epsilon_t)$ is directly observed (see \cite{rt1996}).

For the estimation of SARMA and multivariate SARMA models, the commonly used
estimation method is the quasi-maximum likelihood estimation, which can be also viewed as a nonlinear least squares estimation (LSE).
Given a realization ${X}_1,{X}_2,\dots,{X}_n$ satisfying (\ref{ARMA}), the variable $\epsilon_t(\theta)$ can be approximated, for $0<t\leq n,$ by ${e}_t(\theta)$  defined recursively by
\begin{align}\label{etilde}
{e}_t(\theta)& ={X}_t-\sum_{i=1}^pa_i{X}_{t-i}-\sum_{j=1}^P\mathbf{a}_j{X}_{t-sj}+\sum_{i=1}^p\sum_{j=1}^Pa_i\mathbf{a}_j{X}_{t-sj-i} +
\sum_{i=1}^qb_i{e}_{t-i}(\theta)\nonumber \\&+
\sum_{j=1}^Q\mathbf{b}_j{e}_{t-sj}(\theta)-
\sum_{i=1}^q\sum_{j=1}^Qb_i\mathbf{b}_j{e}_{t-sj-i}(\theta),
\end{align}
where the unknown initial values are set to zero:
${e}_{0}(\theta)=\dots={e}_{1-q-sQ}(\theta)={X}_0=\dots={X}_{1-p-sP}=0$.
The Gaussian quasi-likelihood is given by
\begin{eqnarray*}
{\mathrm{L}}_n(\theta,\sigma^2)=\prod_{t=1}^n\frac{1}{(2\pi)^{1/2}\sqrt{\sigma^2}}\exp\left\{-\frac{{e}^2_t(\theta)}{2\sigma^2}
\right\}.
\end{eqnarray*}
A QMLE of $(\theta,\sigma^2)$ is a measurable solution
$(\hat{\theta}_n,\hat\sigma^2)$ of
\begin{eqnarray*}(\hat{\theta}_n,\hat\sigma^2)
=\arg\min_{\theta,\sigma^2}\left\{
\log(\sigma^2)+\frac{1}{2\sigma^2}Q_n(\theta)\right\}\text{ where }
Q_n(\theta)=\frac{1}{n}\sum_{t=1}^n{e}^2_t(\theta).
\end{eqnarray*}
In all the sequel, we denote  by $\xrightarrow[]{\mathrm{d}}$,
the convergence in distribution.
The symbol $\mathrm{o}_{\mathbb P}(1)$ is used for a sequence of random variables that converges to zero in probability.
Under  the above assumptions, \cite{BMIA2018} showed that $\hat{\theta}_n\to\theta_0\;a.s.$ as $n\to\infty$  and
\begin{eqnarray}\label{NAsymp}
\sqrt{n}\left(\hat{\theta}_n-\theta_0\right)=-J ^{-1}\frac{1}{\sqrt{n}}\sum_{t=1}^n \Upsilon_t+\mathrm{o}_\mathbb P(1)\xrightarrow[n\to\infty]{\mathrm{d}} {\cal
N}(0,\Sigma=J^{-1}IJ^{-1}),
\end{eqnarray}
where
\begin{align*}
J=J(\theta_0)= \frac{2}{\sigma_0^2} \mathbb E \left [ \frac{\partial \epsilon_t(\theta_0)}{\partial\theta}\frac{\partial \epsilon_t(\theta_0)}{\partial\theta'}  \right ],\quad
I=I(\theta_0)=\sum_{h=-\infty}^{+\infty}\cov({\Upsilon}_t,{\Upsilon}_{t-h})\ \text{ and } {\Upsilon}_t=\frac{2}{\sigma_0^2}
 \epsilon_t(\theta_0)\frac{\partial
\epsilon_t(\theta_0)}{\partial\theta}.
\end{align*}
Note that, the existence of the matrix $I(\theta_0)$ is a consequence
of {\bf (A3)} and of Davydov's inequality \cite{davy}.
\section{Diagnostic checking in weak SARMA models} \label{result}
In order to check the validity of the SARMA$(p,q)(P,Q)_s$ model, it is a common practice to
examine the QMLE residuals $\hat{\epsilon}_t=\hat e_t={e}_t(\hat\theta_n)$ where ${e}_t(\theta)$
is given by \eqref{etilde} for all $\theta\in\mathbb{R}^{k_0}$.
For a fixed integer $m\geq1$, consider the vector of residual autocovariances
$${\hat\gamma}_m
=\left(\hat\gamma(1),\dots,\hat\gamma(m)\right)' \  \text{where}\ \; \hat\gamma(h)=\frac{1}{n}\sum_{t=h+1}^{n}\hat e_t\, \hat e_{t-h}
 \  \text{for}\ \;0\leq h<n.$$
In the sequel, we will also need the vector of the first $m$  sample autocorrelations
$$\hat\rho_m =\left(\hat \rho
(1),\dots,\hat \rho(m)\right)' \  \text{where}\ \; \hat \rho(h)=\hat \gamma(h)/\hat \gamma(0).$$
The
statistics \eqref{bp} are usually used to test the following null  hypothesis
\begin{itemize}
\item[($\mathbf{H0}$) : ]  $({X}_t)_{t\in \mathbb Z}$ satisfies a SARMA$(p,q)(P,Q)_s$ representation;
\end{itemize}
against the alternative
\begin{itemize}
\item[($\mathbf{H1}$) :] $({X}_t)_{t\in \mathbb Z}$ does not admit a SARMA$(p,q)(P,Q)_s$ representation or $({X}_t)_{t\in \mathbb Z}$ satisfies a SARMA$(p',q')(P',Q')_s$ representation with $p'>p$ or $q'>q$ or $P'>P$ or $Q'>Q$.
\end{itemize}
\subsection{Asymptotic distribution of the residual autocorrelations}
\label{diagnostic1}
First note that the mixing assumptions {\bf (A3)} entail the asymptotic normality of the "empirical" autocovariances
$$\gamma_m =\left(\gamma(1),\dots,\gamma(m)\right)'
\text{ where } \gamma(h)=\frac{1}{n}\sum_{t=h+1}^{n}\epsilon_t\, \epsilon_{t-h} \  \text{for}\ \;0\leq h<n.$$
It should be noted that $\gamma(h)$ is not a computable statistic because it depends on
the unobserved innovations $\epsilon_t=\epsilon_t(\theta_0)$ except when $p=q=P=Q=0$.
Define the matrix
\begin{eqnarray}\label{xi}
\Xi=\left(\begin{array}{cc} \Sigma&
\Sigma_{\hat\theta_n,\gamma_m}
\\ \Sigma'_{\hat\theta_n,\gamma_m} &
\Sigma_{\gamma_m}
\end{array}\right)=\sum_{h=-\infty}^{+\infty}\mathbb{E}w_tw'_{t-h},\text{ where } w_t=\left(\begin{array}{c}w_{1t}\\w_{2t}\end{array}\right)\in\mathbb R^{k_0+m}
\end{eqnarray}
with $w_{1t}=-J^{-1}\Upsilon_t=-2\sigma_{0}^{-2}\epsilon_tJ^{-1}\left(\partial\epsilon_t/\partial\theta\right)$ and $w_{2t}=\left(\epsilon_{t-1},\dots,\epsilon_{t-m}\right)'\epsilon_t$.
Note that, the existence of the matrix $\Xi$ is a consequence
of {\bf (A3)} and of Davydov's inequality \cite{davy}.

The asymptotic distribution of $\sqrt{n}\hat{\rho}_m$ will be obtained from the
joint asymptotic distribution of $$\sqrt{n}\left(\hat{\theta}'_n-\theta'_0,\gamma'_m\right)'=\frac{1}{\sqrt{n}}\sum_{t=1}^n w_t+\mathrm{o}_{\mathbb{P}}(1)\xrightarrow[n\to\infty]{\mathrm{d}} {\cal
N}(0,\Xi),$$
by applying the central limit theorem for mixing processes (see \cite{herr}).

Now, considering $\hat\gamma(h)$  and $\gamma(h)$ as values of the same
function at the points $\hat\theta_n$ and $\theta_0$, a Taylor
expansion about $\theta_0$ gives
\begin{eqnarray*}
\hat\gamma(h)&=&\gamma(h)+\frac{1}{n}
\sum_{t=h+1}^{n}\left\{\epsilon_{t-h}(\theta)\frac{\partial
\epsilon_{t}(\theta)}{\partial\theta'} + \frac{\partial
\epsilon_{t-h}(\theta)}{\partial\theta'}\epsilon_{t}(\theta)\right\}
_{\theta=\theta_n^*} (\hat\theta_n-\theta_0)+\mathrm{O}_\mathbb{P}(1/n)
\\&=&\gamma(h)+\mathbb{E}\left(\epsilon_{t-h}(\theta_0)\frac{\partial
\epsilon_{t}(\theta_0)}{\partial\theta'}\right)
(\hat\theta_n-\theta_0)+\mathrm{O}_\mathbb{P}(1/n),
\end{eqnarray*}
where $\theta_n^*$ is between $\hat\theta_n$ and $\theta_0.$ The
last equality follows from the consistency of $\hat\theta_n$ and the
fact that $\left(\partial \epsilon_{t-h}/\partial\theta'\right)(\theta_0)$
is not correlated with $\epsilon_t$ when $h\geq0.$ Then for $h=1,\dots,m,$ we have
\begin{eqnarray}\label{taylhatgammapasauto-normalise}
{\hat\gamma}_m :=\left(\hat\gamma(1),\dots,\hat\gamma(m)\right)'=
\gamma_m+\Phi_m(\hat\theta_n-\theta_0)+\mathrm{O}_\mathbb{P}(1/n),
\end{eqnarray}
where
\begin{equation}\label{Phi_m}\Phi_m=\mathbb{E}
\left\{\left(\begin{array}{c}
\epsilon_{t-1}\\\vdots\\\epsilon_{t-m}\end{array}\right)\frac{\partial
\epsilon_{t}}{\partial\theta'}\right\}.
\end{equation}
The following Proposition, which is a generalization of Theorem 4 obtained by \cite{BMIA2018}, gives the limiting distribution of the residual autocovariances and autocorrelations of SARMA models.
\begin{prop}\label{propo1}
When, $p>0$, $P>0$, $q>0$ and $Q>0$, under the above assumptions, we have
$$\sqrt{n}\hat\gamma_m\xrightarrow[n\to\infty]{\mathrm{d}}\mathcal{N}\left(0,
\Sigma_{\hat\gamma_m}\right)\quad \mbox{and}\quad
\sqrt{n}\hat\rho_m\xrightarrow[n\to\infty]{\mathrm{d}}\mathcal{N}\left(0,
\Sigma_{\hat\rho_m}\right)\quad \mbox{where}, $$
\begin{align*}
 \Sigma_{\hat\gamma_m}&=
\Sigma_{\gamma_m}+\Phi_m\Sigma
\Phi'_m+\Phi_m\Sigma_{\hat\theta_n,\gamma_m}
+\Sigma'_{\hat\theta_n,\gamma_m}\Phi'_m\text{ and }
\Sigma_{\hat\rho_m}=\frac{1}{\sigma_0^4}\Sigma_{\hat\gamma_m}.&
\end{align*}
\end{prop}
The proof of this result is similar to that given by  \cite{BMIA2018} for Theorem 4.

The asymptotic variance matrices $\Sigma_{\hat\gamma_m}$ and $\Sigma_{\hat\rho_m}$ depend on the unknown matrices $\Xi$, $\Phi_m$ and  the scalar $\sigma_0^2$.  Matrix  $\Phi_m$ and $\sigma_0^2$ can be estimated by its empirical counterpart, respectively
$$\hat{\Phi}_m=\frac{1}{n}\sum_{t=1}^n
\left\{\left(\hat{\epsilon}_{t-1},\dots,\hat{\epsilon}_{t-m}\right)'\frac{\partial
\hat{\epsilon}_{t}}{\partial\theta'}\right\}\text{ and }\hat{\sigma}^2=\hat\gamma(0)=\frac{1}{n}\sum_{t=1}^n\hat{\epsilon}_t^2.$$
Note that the matrix $(2\pi)^{-1}\Xi$  is the spectral density at frequency zero of the process $(w_t)$, thus an estimator of $\Xi$ is given in Theorem 6 of \cite{BMIA2018}.
Other estimators of such long-run variances are available in the literature (see for instance \cite{a_econ}, \cite{berk}, \cite{haan}, \cite{newey}, for general references).
For the numerical illustrations presented in this paper, we used a Vector AR (VAR) spectral
estimator given in Theorem 6 of \cite{BMIA2018} consisting in: i) fitting VAR$(r)$ models for $r= 1,\dots,r_{\max}$ to the series $\hat{w}_t$, $t=1,\dots,n$,
where $\hat{w}_t$ is obtained by replacing $\theta_0$ by $\hat{\theta}_n$ in
$w_t$; ii) selecting the order $r$ which minimizes an information criterion and approximating $\Xi$ by
$(2\pi)$ times the spectral density at frequency zero of the estimated VAR$(r)$ model. Hereafter, we used
the AIC model selection criterion with $r_{\max}=5$.

From Proposition \ref{propo1} we can deduce the following result, which gives  the exact  limiting distribution of the standard portmanteau statistics \eqref{bp} under general assumptions on the innovation process of the fitted SARMA$(p,q)(P,Q)_s$ model.
\begin{theo}\label{theolimitdistBP}
Under Assumptions in Proposition \ref{propo1} and ${\bf (H0)}$,  the
statistics $Q_m^{\textsc{lb}}$ and $Q_m^{\textsc{bp}}$ converge in distribution, as $n\rightarrow\infty,$ to
$$Z_m(\xi_m)=\sum_{i=1}^{m}\xi_{i,m}Z_i^2$$ where
$\xi_m=(\xi_{1,m},\dots,\xi_{m,m})'$ is the vector of the
eigenvalues of the matrix $\Sigma_{\hat\rho_m}=\sigma_0^{-4}\Sigma_{\hat\gamma_m}$  and $Z_1,\dots,Z_{m}$ are independent
$\mathcal{N}(0,1)$ variables.
\end{theo}
As in \cite{BMIA2018}, Theorem \ref{theolimitdistBP} shows that for the asymptotic distribution of $Q_m^{\textsc{lb}}$ and $Q_m^{\textsc{bp}}$, the
$\chi_{m-k_0}^2$  approximation is no longer valid in the framework of weak SARMA$(p,q)(P,Q)_s$ models. The true asymptotic distribution depends on nuisance parameters involving $\sigma_0^2$, the matrix $\Phi_m$ and the elements of $\Xi$. Consequently, in order to obtain the asymptotic distribution of the portmanteau statistics \eqref{bp} under weak assumptions on the noise, one needs a consistent estimator of the asymptotic covariance matrix $\Sigma_{\hat\rho_m}$.
We let $\hat\Sigma_{\hat\rho_m}$ the matrix obtained by replacing $\Xi$ by $\hat\Xi$, $\Phi_m$ by $\hat\Phi_m$ and $\sigma_0^2$ by $\hat\sigma^2$ in $\Sigma_{\hat\rho_m}$. Denote by $\hat\xi_m=(\hat\xi_{1,m},\dots,\hat\xi_{m,m})'$ the
vector of the eigenvalues of $\hat\Sigma_{\hat\rho_m}$. At the asymptotic level $\alpha$, the LB (Ljung-Box) test (resp. the BP (Box-Pierce) test) consists in rejecting the adequacy of the weak SARMA$(p,q)(P,Q)_s$ model when $${Q}_m^{\textsc{lb}}>S_m(1-\alpha)\quad (\text{ resp. }{Q}_m^{\textsc{bp}}>S_m(1-\alpha)),$$
where $S_m(1-\alpha)$ is such that
$\mathbb{P}\left\{Z_m(\hat\xi_m)>S_m(1-\alpha)\right\}=\alpha$.
We emphasize the fact that the proposed modified versions of the Box-Pierce and Ljung-Box statistics are more difficult to implement because their critical values have to be computed from the data.
\subsection{Self-normalized asymptotic distribution of the residual autocorrelations}
\label{diagnostic2}
The nonparametric kernel estimator (see \cite{a_econ,newey}), used to estimate the matrix $\Xi$  causes serious difficulties regarding the choice of the sequence of weights.
The parametric approach based on an autoregressive  estimate of the spectral density of $w_t$  studied for instance by \cite{berk,yac,yac2,BMIA2018,haan}  is also facing the problem of choosing the truncation  parameter. So the choice of the order of truncation is often crucial and difficult.
In this section, we  propose  as in \cite{BMS2018} an alternative method where we do not estimate an asymptotic covariance matrix. It is based on a self-normalization
based approach to construct a test-statistic which is asymptotically distribution-free under the null hypothesis (see \cite{BMS2018}, for a reference in the ARMA cases).
The idea comes from \cite{lobato} and has been already extended by \cite{kl2006,s2010JRSSBa,s2010JRSSBb,shaox} to more general frameworks.
See also \cite{s2016} for a review on some recent developments on the inference of time series data using the self-normalized approach.
In this  case, the critical values are not computed from the data since they are tabulated. In some sense,  this method is finally closer to the standard method in which the critical values are simply deduced from a $\chi^2$-table.

%
We denote $\Lambda$ the matrix in $\mathbb R^{m\times(k_0+m)}$ defined in block
formed by
${\Lambda}=({\Phi}_m {\vert} I_{m})$, where $I_m$ is the  identity matrix of order $m$.
In view of \eqref{NAsymp} and \eqref{taylhatgammapasauto-normalise}, we  deduce that
\begin{align}\label{tayl}
\sqrt{n}\ \hat\gamma_m & = \frac{1}{\sqrt{n}}\sum_{t=1}^n\ \Lambda w_t + \mathrm{o}_{\mathbb P}(1).
\end{align}
Contrarily to Subsection \ref{diagnostic1}, we do not rely on the classical method that would consist in estimating the asymptotic covariance matrix of $\Lambda w_t$. We need to apply the functional central limit theorem holds for the process $w=(w_t)_{t\ge 1}$ (see Lemma 1 in \cite{lobato}).

Finally, we define the normalization matrix $C_{m}\in\mathbb R^{m\times m}$  by
\begin{equation*}
C_{m} =\frac{1}{n^2}\sum_{t=1}^{n}{S}_{t}{S}'_{t} \text{ where }
{S}_t =\sum_{j=1}^{t}\left({\Lambda}{w}_j-\Lambda\bar{w}\right) \text{ with }\bar{w}=\frac{1}{n}\sum_{t=1}^nw_t.
\end{equation*}
To ensure the invertibility of the normalization matrix $C_{m}$ which is proved  in Lemma 6 of \cite{BMS2018}, we
need the following technical assumption on the distribution of $\epsilon_t$.
\begin{itemize}
\item[\hspace*{1em} {\bf (A4):}]
\hspace*{1em} {The process $(\epsilon_t)_{t\in\mathbb Z}$ has a positive density on some neighborhood of zero}.
\end{itemize}
 Let $(B_K(r))_{r\ge 0}$ be a $K$-dimensional Brownian motion starting from $0$.  For $K\ge 1$, we denote  $\mathcal U_K$ the random variable defined by
\begin{equation}\label{UkVk}
\mathcal U_K={B}'_{K}(1){V}_{K}^{-1}{B}_{K}(1)\text{ where }
{V}_{K}=\int_0^1\left ( {B}_{K}(r)-r{B}_{K}(1)\right)\left({B}_{K}(r)-r{B}_{K}(1)\right)'dr.
\end{equation}
The following theorem states the asymptotic distributions of the sample autocovariances and autocorrelations.
\begin{theo}\label{sn1}
We assume that $p>0$, $q>0$, $P>0$ or $Q>0$. Under Assumptions of Proposition \ref{propo1}, {\bf (A4)} and  under the null hypothesis {\bf (H0)}, we have
\begin{align*}
n\, \hat\gamma_m'C_{m}^{-1}\hat\gamma_m \xrightarrow[n\to\infty]{\mathrm{d}} \mathcal{U}_{m}
\text{ and }n\,\sigma_0^4\hat{{\rho}}_m'C_{m}^{-1}\hat{{\rho}}_m
\xrightarrow[n\to\infty]{\mathrm{d}} \mathcal{U}_{m}
\end{align*}
\end{theo}
The proof of this result is postponed to Section \ref{proof}.

Of course, the above theorem is useless for practical purpose,  because it does not involve any observable quantities. 
In practice, one has to replace the matrix $C_{m}$ and the variance of
the noise $\sigma_{0}^2$  by their empirical or observable counterparts.
The matrix
$J$ can be easily estimated by his empirical counterpart
\begin{equation*}
\hat{{J}}=\frac{2}{\hat\sigma^2}\frac{1}{n} \sum_{t=1}^n \frac{\partial
{e}_t(\hat\theta_n)}{\partial
\theta}\frac{\partial
{e}_t(\hat\theta_n)}{\partial \theta'}.
\end{equation*}
Thus we define
\begin{equation*}
 \hat\Lambda =  \Big ( \hat \Phi_m {\vert} I_m  \Big ) \text{ and }
  \hat w_t = \left (-2 \hat J_n^{-1}\frac{\partial e_t(\hat\theta_n) }{\partial \theta'}\frac{1}{\hat\sigma^2}\hat e_t,\hat e_t\hat e_{t-1},\dots,\hat e_t\hat e_{t-m}\right )'.
\end{equation*}
Finally we denote the normalization matrix $\hat C_m\in\mathbb R^{m\times m}$ by
\begin{equation*}
\hat C_{m} =\frac{1}{n^2}\sum_{t=1}^{n}\hat S_{t}\hat S'_{t}  \text{  where  }
\hat S_t =\sum_{j=1}^{t}\left(\hat\Lambda\hat w_j-\hat\Lambda\bar{\hat{w}}\right) \text{ with }\bar{\hat{w}}=\frac{1}{n}\sum_{t=1}^n\hat{w}_t.
\end{equation*}
The following result is the applicable counterpart of Theorem \ref{sn1}.
\begin{theo}\label{sn2}
Assume that $p>0$, $q>0$, $P>0$ or $Q>0$. Under Assumptions of Theorem \ref{sn1}, we have
\begin{align*}
n\, \hat\gamma_m'\hat C_m^{-1}\hat\gamma_m \xrightarrow[n\to\infty]{\mathrm{d}} \mathcal{U}_m\text{ and }
Q_m^\textsc{sn}=n\,\hat\sigma^4 \,\hat\rho_m'\hat C_m^{-1}\hat\rho_m  \xrightarrow[n\to\infty]{\mathrm{d}} \mathcal{U}_m.
\end{align*}
\end{theo}
The proof of this result is postponed to Section \ref{proof}.

Based on the above result, we propose a modified version of the Ljung-Box statistic when one uses the statistic
\begin{equation}\label{snLB}
\tilde Q_m^\textsc{sn}= n\, \hat\sigma^4  \,\hat\rho_m'  D^{1/2}_{n,m} \hat C_m^{-1}D^{1/2}_{n,m}\hat\rho_m,
\end{equation}
where the matrix $D_{n,m}\in\mathbb R^{m\times m}$ is diagonal with $((n+2)/(n-1),...,(n+2)/(n-m))$ as diagonal terms.
\section{Numerical illustrations}\label{ne}
In this section, by means of Monte Carlo experiments, we investigate the finite sample
properties of the modified version of the portmanteau tests that we introduced in this work.
The numerical illustrations of this section are made  with the open source
statistical software R (see R Development Core Team, 2017) or (see http://cran.r-project.org/).
\subsection{Simulated models}
\label{univARMA}
First of all, we introduce the models that we simulate and we indicate the conventions that we adopt in the discussion and in the tables:
\begin{itemize}
\item $\mathrm{{LB}_{\textsc{w}}}$ and $\mathrm{BP_{\textsc{w}}}$ refer to  modified LB and BP tests using  ${Q}_m^\textsc{lb}$ and $Q_m^\textsc{bp}$  in Section \ref{diagnostic1}
\item $\mathrm{{LB}_{\textsc{s}}}$ and $\mathrm{BP_{\textsc{s}}}$ refer to  LB and BP tests using the standard statistics \eqref{bp}.
\item $\mathrm{{{LB}}_{\textsc{sn}}}$ and  $\mathrm{{{BP}}_{\textsc{sn}}}$ refer to modified tests using the self-normalized statistics in Section \ref{diagnostic2}
\end{itemize}
To generate the strong and the weak SARMA models, we consider the
following SARMA$(0,1)(0,1)_s$ model
\begin{eqnarray}
\label{ARMA01MonteCarlo}
X_{t}=\epsilon_{t}-b_{01}\epsilon_{t-1}-\mathbf{b}_{01}\epsilon_{t-s}+b_{01}\mathbf{b}_{01}\epsilon_{t-s-1},
\end{eqnarray}
with $\theta_0=(b_{01},\mathbf{b}_{01})'=(-0.6,-0.7)'$ and the innovation process  $(\epsilon_t)$ follows a strong or weak white noise.

The  generalized autoregressive conditional heteroscedastic (GARCH) models is an important example of weak white noises in the univariate case (see \cite{FZ2010}).
So we first assume that in (\ref{ARMA01MonteCarlo}) the
innovation process $\epsilon$ is the following ARCH$(1)$ model
defined by
\begin{equation} \label{bruitARCH}
\left\{\begin{array}{l}\epsilon_{t}=\sigma_t\eta_{t}\\
\sigma_t^2=1+\alpha_1\epsilon_{t-1}^2 
\end{array}\right.
\end{equation}
where  $(\eta_t)_{t\ge 1}$ is a sequence of iid standard Gaussian random variables. 
To generate the strong SARMA, we assume that in
(\ref{ARMA01MonteCarlo}) the innovation process follows (\ref{bruitARCH}) with $\alpha_1=0$.
\subsection{Empirical size}
We first simulate  $N=1,000$ independent trajectories of size $n=2,000$ of models (\ref{ARMA01MonteCarlo}). The same series is  partitioned as two series of sizes $n=500$ and $n=2,000$.
For each of these $N$ replications,  we use the quasi-maximum likelihood estimation method to estimate the coefficient $\theta_0$ and we apply portmanteau tests to the residuals for different values of $m\in\{4,8,12,15,18,20\}$, where $m$ is the number of autocorrelations used in the portmanteau test statistic.
For the nominal level $\alpha=5\%$, the
empirical  size over the $N$ independent replications should
vary between the significant limits 3.6\% and 6.4\% with probability
95\% and belong to $[3.2\%,6.9\%]$ with a probability 99\%.
When   the relative rejection frequencies  are outside the 95\%
significant limits, they are displayed in bold type and they are underlined when they  are outside the 99\%
significant limits in Tables \ref{tab1ARMA} and \ref{tab2ARMA}.

For the standard  Box-Pierce test, the model is therefore rejected when the statistic $Q_m^{\textsc{bp}}$ or $Q_m^{\textsc{lb}}$ is larger than $\chi_{(m-2)}^2(0.95)$  in a
SARMA$(0,1)(0,1)_s$ case (see \cite{McLeod1978}). Consequently the empirical size is not available (n.a.) for the statistic $Q_m^{\textsc{bp}}$ or $Q_m^{\textsc{lb}}$
because they are not applicable for $m\leq 2$.
For the proposed self-normalized test $\mathrm{BP}_\textsc{sn}$ or $\mathrm{LB}_\textsc{sn}$, the model is rejected when the statistic $
Q_m^\textsc{sn}$ or $\tilde{Q}_m^\textsc{sn}$ is larger than $\mathcal{U}_m(0.95)$, where the critical values $\mathcal{U}_K(0.95)$ (for $K=1,\dots,20$) are tabulated  in Lobato (see Table 1 in \cite{lobato}).

Table \ref{tab1ARMA}  displays the relative rejection frequencies of the null hypothesis ${\bf (H0)}$ that the  data generating process (DGP for short) follows a
strong SARMA model (\ref{ARMA01MonteCarlo})--(\ref{bruitARCH}) with $\alpha_1=0$, over the $N$ independent replications.
When the seasonal period is $s=4$, for all tests, the percentages of rejection belong to the confident interval with probabilities 95\%
and 99\%, except for $\mathrm{{LB}}_\textsc{s}$ and ${\mathrm{BP}}_\textsc{s}$ when $m=4$. Consequently all these tests well control the error of first kind. In contrast,
when $s=12$, our proposed tests well also control the error of first kind (except for $\mathrm{{LB}}_\textsc{w}$ and ${\mathrm{BP}}_\textsc{w}$ when $n=500$)
contrarily to the standard tests $\mathrm{{LB}}_\textsc{s}$, ${\mathrm{BP}}_\textsc{s}$ for all sizes.
We draw the conclusion that, in this strong SARMA case, the proposed modified version may be clearly preferable to the standard ones.

\begin{figure}[hbt]
\vspace*{12cm} \protect \includegraphics{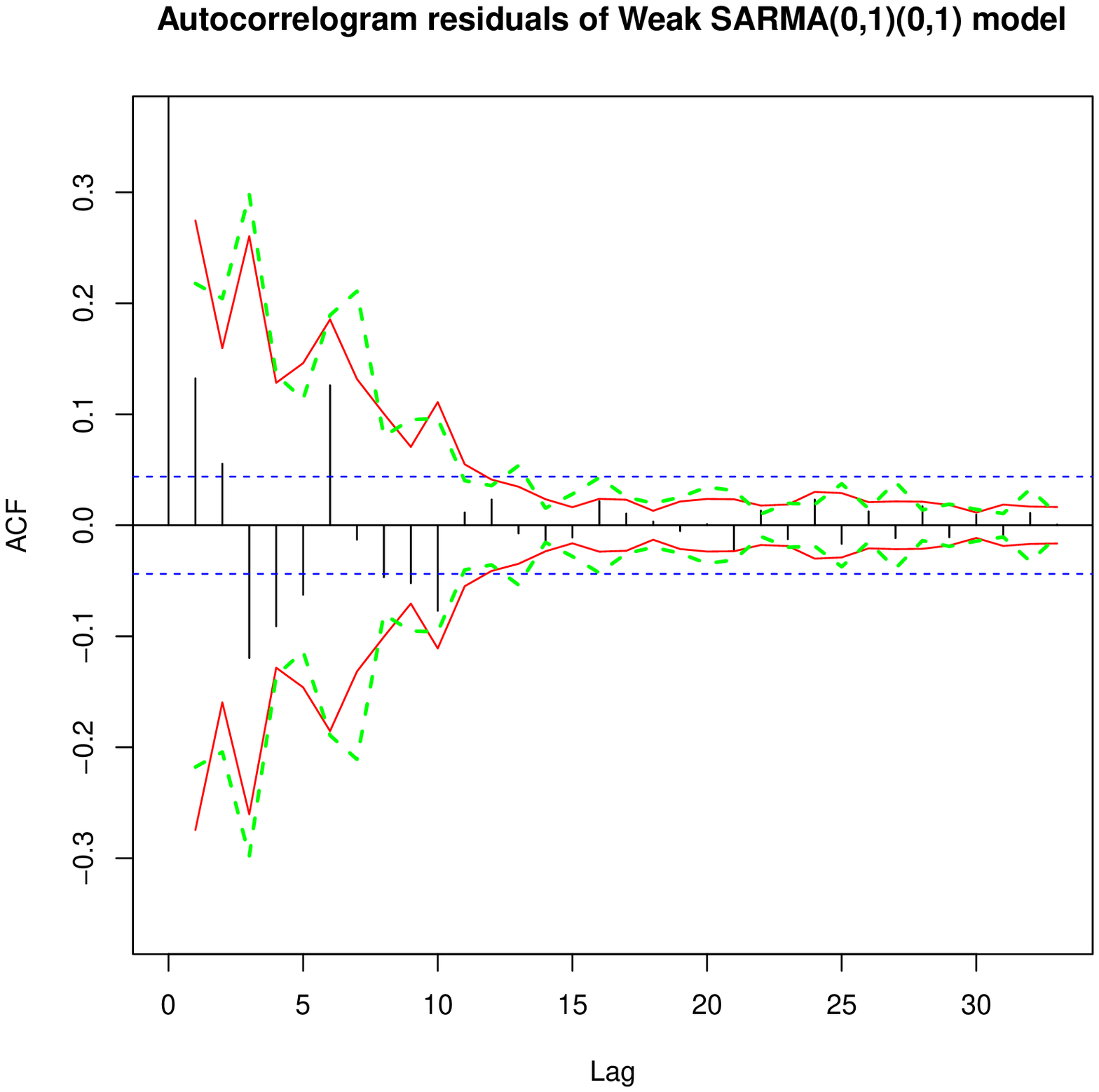} \vspace*{2.5cm}
\caption{\label{fig} {\footnotesize Autocorrelation  of a realization of size $n=2,000$ for weak SARMA$(0,1)(0,1)_{12}$ model
(\ref{ARMA01MonteCarlo})-(\ref{bruitARCH}) with $\alpha_1=1.3$.
The horizontal dotted lines (blue color) correspond to the 5\% significant limits obtained under the strong SARMA assumption.
The solid lines (red color) and  dashed lines (green color) correspond also  to the 5\% significant limits under the weak SARMA assumption.
The full lines  correspond to the asymptotic significance limits for the residual autocorrelations
obtained in Proposition~\ref{propo1}. The dashed lines (green color) correspond to the self-normalized asymptotic significance limits for the residual autocorrelations
obtained in Theorem~\ref{sn2}.}}
\end{figure}

Now, we repeat the same experiments on a weak SARMA models.
As expected, Table \ref{tab2ARMA}  shows that the
standard  $\mathrm{{LB}}_\textsc{s}$ or $\mathrm{{BP}}_\textsc{s}$ test poorly performs in assessing the adequacy of this particular weak SARMA model.
It can be seen that: 1) the observed relative rejection frequencies of $\mathrm{{LB}}_\textsc{s}$ and $\mathrm{{BP}}_\textsc{s}$ are
definitely outside the significant limits, 2) the errors of the first kind are only globally well controlled by the proposed tests, for all $s$ when $n$ is large.
We also tried the case where the ARCH$(1)$ model (\ref{bruitARCH}) have infinite fourth moments. As showing in Figure~\ref{fig}, the results are qualitatively
similar to what we observe here.

Figure~\ref{fig} displays the residual autocorrelations  of a realization of size $n=2,000$ for weak SARMA$(0,1)(0,1)_{12}$ model
(\ref{ARMA01MonteCarlo})-(\ref{bruitARCH}) with $\alpha_1=1.3$ and their 5\% significance limits under the strong  SARMA$(0,1)(0,1)_{12}$
and weak  SARMA$(0,1)(0,1)_{12}$ assumptions. This figure confirms clearly the conclusions drawn from Table  \ref{tab2ARMA}.
The horizontal dotted lines (blue color) correspond to the 5\% significant limits obtained under the strong SARMA assumption.
The solid lines (red color) and  dashed lines (green color) correspond also  to the 5\% significant limits under the weak SARMA assumption.
The full lines  correspond to the asymptotic significance limits for the residual autocorrelations
obtained in Proposition~\ref{propo1}. The dashed lines (green color) correspond to the self-normalized asymptotic significance limits for the residual autocorrelations
obtained in Theorem~\ref{sn2}.

In these Monte Carlo experiments, we illustrate that the proposed test statistics  have reasonable finite sample performance. Under nonindependent errors,
it appears that the standard test statistics are generally non reliable, overrejecting  severely, while the proposed tests statistics offer satisfactory levels.
Even for independent errors, they seem  preferable to the standard ones, when the number $m$ of autocorrelations is small and when $s=12$. Moreover, the error of first kind is well controlled.
Contrarily to the standard tests  based on $\mathrm{{BP}}_{\textsc{s}}$ or $\mathrm{{{LB}}_{\textsc{s}}}$, the proposed tests can be used safely for $m$ small (see for instance Figure~\ref{fig}).
For all these above reasons, we think that the modified
versions that we propose in this paper are preferable to the standard ones for diagnosing SARMA models under nonindependent errors.
\begin{table}[hbt]
 \caption{\small{Empirical size (in \%) of the modified and standard versions
 of the LB and BP tests in the case of SARMA$(0,1)(0,1)_s$ model
(\ref{ARMA01MonteCarlo})-(\ref{bruitARCH}) with $\alpha_1=0.$
The nominal asymptotic level of the tests is $\alpha=5\%$.
The number of replications is $N=1,000$. }}
\begin{center}
\begin{tabular}{c c ccc ccc c}
\hline\hline \\
$s$& Length $n$ & Lag $m$ & $\mathrm{{LB}}_{\textsc{sn}}$&$\mathrm{BP}_{\textsc{sn}}$&$\mathrm{{LB}}_{\textsc{w}}$&$\mathrm{BP}_{\textsc{w}}$&$\mathrm{{LB}}_{\textsc{s}}$&$\mathrm{BP}_{\textsc{s}}$
\vspace*{0.2cm}\\\hline
&& $4$&4.3& 4.0& {\bf 6.8}& {\bf 6.8}& \underline{{\bf 8.7}}&\underline{{\bf 8.7}}\\
&& $8$&5.8& 5.4& 5.8& 5.6 &\underline{{\bf 7.2}} &{\bf 6.7}\\
4 &$n=500$& $12$&4.9& 4.8& 6.0& 5.4& \underline{{\bf 7.2}}& 6.3\\
 &&$15$&6.3& 5.3& 5.1& 4.7& 5.8& 5.3\\
 &&$18$&5.7& 5.1& 5.4& 4.3 &5.7 &4.9\\
 &&$20$&5.7& 5.0& 5.2& 3.8& 5.4& 4.1\\
  \cline{2-9}
&& $4$&4.1& 4.1& 5.7& 5.6&\underline{{\bf 7.5}}&\underline{{\bf 7.5}}\\
&& $8$&4.2 &4.2& 4.7& 4.7& 5.9& 5.6\\
4 &$n=2,000$& $12$&4.9& 4.8& 3.9& 3.9& 4.6& 4.4\\
 &&$15$&4.6& 4.6& 4.9& 4.4& 4.9& 4.5\\
 &&$18$&4.0& 3.8& 4.4& 4.2& 4.6& 4.1\\
 &&$20$&4.5& 4.5& 5.1& 5.0& 5.4& 5.1\\

\hline
&& $4$&4.7 & 4.5& \underline{{\bf 12.8}}&\underline{{\bf 12.1}} &\underline{{\bf 20.8}} &\underline{{\bf 20.6}}\\
&& $8$& {\bf 6.5} & 6.4 &\underline{{\bf 12.1}}&\underline{{\bf 11.4}}&\underline{{\bf 15.1}} &\underline{{\bf 14.6}}\\
12 &$n=500$& $12$&5.4 & 5.1&\underline{{\bf 10.7}}&\underline{{\bf 10.1}}&\underline{{\bf 11.7}}&
\underline{{\bf 10.8}}\\
 &&$15$&6.0& 5.6&\underline{{\bf 8.7}}&\underline{{\bf 8.4}}&\underline{{\bf 9.8}}&\underline{{\bf 9.4}}\\
 &&$18$&5.1& 4.5& \underline{{\bf 8.3}}&\underline{{\bf 7.5}}&\underline{{\bf 9.7}}&\underline{{\bf 8.3}}\\
 &&$20$&5.5& 4.3&\underline{{\bf 8.6}}&\underline{{\bf 7.6}}&\underline{{\bf 9.8}}&\underline{{\bf 8.4}}\\
 \cline{2-9}
&& $4$&4.4 & 4.3 &\underline{{\bf 7.2}} &\underline{{\bf 7.0}} &\underline{{\bf 14.0}}&\underline{{\bf 13.7}}\\
&& $8$&4.9& 4.9& 5.8& 5.7&\underline{{\bf 9.9}}&\underline{{\bf 9.8}}\\
12 &$n=2,000$& $12$&4.5& 4.4& 6.0& 5.6&\underline{{\bf 7.0}}&{\bf 6.6}\\
 &&$15$&4.5& 4.5& 6.1& 6.0&\underline{{\bf 7.0}}&{\bf 6.9}\\
 &&$18$&5.4& 5.3& 6.0& 5.7&{\bf 6.9}&{\bf 6.5}\\
 &&$20$&4.9& 4.8& 6.3& 6.1 &\underline{{\bf 7.5}}&\underline{{\bf 7.1}}\\
\hline\hline
\\
\end{tabular}
\end{center}
\label{tab1ARMA}
\end{table}

\begin{table}[hbt]
 \caption{\small{Empirical size (in \%) of the modified and standard versions
 of the LB and BP tests in the case of SARMA$(0,1)(0,1)_s$ model
(\ref{ARMA01MonteCarlo})-(\ref{bruitARCH}) with $\alpha_1=0.45.$
The nominal asymptotic level of the tests is $\alpha=5\%$.
The number of replications is $N=1,000$. }}
\begin{center}
\begin{tabular}{c c ccc ccc c}
\hline\hline \\
$s$& Length $n$ & Lag $m$ & $\mathrm{{LB}}_{\textsc{sn}}$&$\mathrm{BP}_{\textsc{sn}}$&$\mathrm{{LB}}_{\textsc{w}}$&$\mathrm{BP}_{\textsc{w}}$&$\mathrm{{LB}}_{\textsc{s}}$&$\mathrm{BP}_{\textsc{s}}$
\vspace*{0.2cm}\\\hline
&& $4$&\underline{{\bf 3.1}} &\underline{{\bf 2.9}} & 4.5 & 4.3&\underline{{\bf 18.3}} &\underline{{\bf 18.3}}\\
&& $8$&\underline{{\bf 2.8}} &\underline{{\bf 2.6}} &{\bf 3.3}&{\bf 3.2}&\underline{{\bf 11.7}} &\underline{{\bf 11.3}}\\
4 &$n=500$& $12$&{\bf 3.5}&{\bf 3.2}&\underline{{\bf  2.8}}&\underline{{\bf 2.4}}&\underline{{\bf  9.7}} &\underline{{\bf 9.3}}\\
 &&$15$&\underline{{\bf 2.8}} &\underline{{\bf  2.3}} &\underline{{\bf 2.9}} &\underline{{\bf 2.7}} &\underline{{\bf 10.0}} &\underline{{\bf 9.0}}\\
 &&$18$&\underline{{\bf 1.9}}&\underline{{\bf  1.6}}& \underline{{\bf 2.7}}&\underline{{\bf  2.4}}&\underline{{\bf 9.1}}&\underline{{\bf 7.8}}\\
 &&$20$&\underline{{\bf 1.7}}&\underline{{\bf 1.5}}&\underline{{\bf 2.6}}&\underline{{\bf 2.3}}&\underline{{\bf 9.0}}&\underline{{\bf 7.9}}\\
  \cline{2-9}
&& $4$&5.5 & 5.4 & 4.2 & 4.2&\underline{{\bf 19.9}}&\underline{{\bf 19.9}}\\
&& $8$&4.8 & 4.8 & 3.6 & 3.6& \underline{{\bf 13.9}}&\underline{{\bf 13.7}}\\
4 &$n=2,000$& $12$&4.5 & 4.3 & 3.8 & 3.6&\underline{{\bf 11.4}} &\underline{{\bf 11.3}}\\
 &&$15$&4.8&  4.8 & 3.6 & 3.6&\underline{{\bf 10.7}}&\underline{{\bf 10.5}}\\
 &&$18$&3.7 & 3.7 & {\bf 3.4} &{\bf 3.4}&\underline{{\bf 10.2}}&\underline{{\bf 9.8}}\\
 &&$20$&{\bf 3.2}&\underline{{\bf 3.1}}&{\bf 3.5}&{\bf 3.5}&\underline{{\bf  9.4}}&\underline{{\bf  9.1}}\\

\hline
&& $4$&3.6 &{\bf 3.5} & \underline{{\bf 9.1}} &\underline{{\bf 8.9}}&\underline{{\bf 29.6}}&\underline{{\bf 29.1}}\\
&& $8$&4.0 & 3.8 & 6.4 & 6.0 &\underline{{\bf 20.4}} &\underline{{\bf 19.7}}\\
12 &$n=500$& $12$&4.2 &{\bf 3.5} & 5.5 & 5.2&\underline{{\bf 15.4}} &\underline{{\bf 14.6}}\\
 &&$15$&\underline{{\bf 2.2}} &\underline{{\bf  1.9}} & 5.4 & 4.8& \underline{{\bf 14.6}} &\underline{{\bf 13.5}}\\
 &&$18$&\underline{{\bf 1.7}} &\underline{{\bf 1.4}} & 4.3 & 4.2&\underline{{\bf 13.3}}&\underline{{\bf 12.0}}\\
 &&$20$&\underline{{\bf 1.8}}&\underline{{\bf 1.6}}&  4.9 & 4.2&\underline{{\bf 12.9}}&\underline{{\bf  11.8}}\\
 \cline{2-9}
&& $4$&4.4&  4.4 & 4.5 & 4.5&\underline{{\bf 27.1}}&\underline{{\bf  27.1}}\\
&& $8$&4.8 & 4.7 & 4.1&  3.9& \underline{{\bf 20.3}}&\underline{{\bf 20.1}}\\
12 &$n=2,000$& $12$&3.8 & 3.8 & 4.6&  4.4&\underline{{\bf 15.3}}&\underline{{\bf 15.2}}\\
 &&$15$&4.6 & 4.5 &{\bf 3.5} &{\bf 3.4}& \underline{{\bf 13.9}}&\underline{{\bf  13.7}}\\
 &&$18$&3.8 & 3.7 & 3.6& {\bf 3.5}&\underline{{\bf 12.6}}&\underline{{\bf 12.4}}\\
 &&$20$&{\bf 3.5}  &{\bf 3.3} & 4.1 & 3.8&\underline{{\bf 11.9}}&\underline{{\bf 11.8}}\\
\hline\hline
\\
\end{tabular}
\end{center}
\label{tab2ARMA}
\end{table}

\subsection{Empirical power}\label{emp-power}
In this section we repeat the same experiments as in Section \ref{univARMA} to examine the power of the tests for the null hypothesis of a SARMA$(0,1)(0,1)_s$ against the following
SARMA$(1,1)(0,1)_s$ alternative defined by
\begin{eqnarray}
X_{t}&=&a_{01}X_{t-1}+\epsilon_{t}-b_{01}\epsilon_{t-1}-\mathbf{b}_{01}\epsilon_{t-s}+b_{01}\mathbf{b}_{01}\epsilon_{t-s-1},\label{ARMA21MonteCarlo}
\end{eqnarray}
with $\theta_0=(a_{01},b_{01},\mathbf{b}_{01})'=(0.8,-0.6,-0.7)'$ and
where the innovation process  $\epsilon$ follows a strong or weak white noise introduced in Section \ref{univARMA}.
For each of these $N$ replications we fit a SARMA$(0,1)(0,1)_s$ models
and perform standard and modified tests based on $m=4,8,12,15$, $18$ and $20$ residual autocorrelations.

Tables \ref{tabfortpuissance} and \ref{tabGARCHpuissance} compare the empirical powers of Model \eqref{ARMA21MonteCarlo}-\eqref{bruitARCH}
with $\alpha_1=0$ and $\alpha_1=0.45$ respectively over the $N$ independent replications.
For these particular strong and weak SARMA models,  we notice that the standard $\mathrm{{BP}}_{\textsc{s}}$ and
$\mathrm{{{LB}}_{\textsc{s}}}$ and our proposed  tests have very similar powers except for $\mathrm{{BP}}_{\textsc{sn}}$ and $\mathrm{{{LB}}_{\textsc{sn}}}$ when  $n=500$ in the weak case.

\begin{table}[hbt]
\caption{\small{Empirical power (in \%) of the modified and standard versions
 of the LB and BP tests in the case of SARMA$(0,1)(0,1)_s$ model \eqref{ARMA21MonteCarlo}-\eqref{bruitARCH} with $\alpha_1=0$.
The number of replications is $N=1,000$. }}
\begin{center}
\begin{tabular}{c c ccc ccc c}
\hline\hline \\
$s$& Length $n$ & Lag $m$ & $\mathrm{{LB}}_{\textsc{sn}}$&$\mathrm{BP}_{\textsc{sn}}$&$\mathrm{{LB}}_{\textsc{w}}$&$\mathrm{BP}_{\textsc{w}}$&$\mathrm{{LB}}_{\textsc{s}}$&$\mathrm{BP}_{\textsc{s}}$
\vspace*{0.2cm}\\\hline
&& $4$&98.9&  98.9& 100.0 &100.0 &100.0 &100.0\\
&& $8$&98.1 & 98.0& 100.0& 100.0 &100.0& 100.0\\
4 &$n=500$& $12$&96.3 & 96.3& 100.0& 100.0 &100.0 &100.0\\
 &&$15$&94.9 & 94.9 & 99.9 & 99.9& 100.0& 100.0\\
 &&$18$&93.7 & 93.5 & 99.8 & 99.7& 100.0 &100.0\\
 &&$20$&91.6 & 91.4 & 99.7 & 99.7 &100.0& 100.0\\
  \cline{2-9}
&& $4$&100.0 &100.0& 100.0& 100.0 &100.0& 100.0\\
&& $8$&100.0 &100.0& 100.0& 100.0 &100.0& 100.0\\
4 &$n=2,000$& $12$&100.0 &100.0& 100.0& 100.0 &100.0& 100.0\\
 &&$15$&100.0 &100.0& 100.0& 100.0 &100.0& 100.0\\
 &&$18$&100.0 &100.0& 100.0& 100.0 &100.0& 100.0\\
 &&$20$&100.0 &100.0& 100.0& 100.0 &100.0& 100.0\\

\hline
&& $4$&98.1 & 98.1& 100.0& 100.0& 100.0 &100.0\\
&& $8$&97.7&  97.7 & 99.9 & 99.9& 100.0& 100.0\\
12 &$n=500$& $12$&96.6 & 96.6 & 99.9 & 99.9& 100.0 &100.0\\
 &&$15$&95.5 & 95.4 & 99.9 & 99.9& 100.0& 100.0\\
 &&$18$&92.9 & 92.7 & 99.9 & 99.9& 100.0& 100.0\\
 &&$20$&90.9 & 90.5 & 99.9 & 99.9& 100.0 &100.0\\
 \cline{2-9}
&& $4$&100.0& 100.0& 100.0& 100.0& 100.0& 100.0\\
&& $8$&100.0& 100.0& 100.0& 100.0 &100.0& 100.0\\
12 &$n=2,000$& $12$&100.0 &100.0& 100.0& 100.0 &100.0& 100.0\\
 &&$15$&100.0& 100.0& 100.0& 100.0& 100.0& 100.0\\
 &&$18$&100.0& 100.0& 100.0& 100.0& 100.0& 100.0\\
 &&$20$&100.0 &100.0& 100.0& 100.0 &100.0& 100.0\\
\hline\hline
\\
\end{tabular}
\end{center}
\label{tabfortpuissance}
\end{table}

\begin{table}[hbt]
\caption{\small{Empirical power (in \%) of the modified and standard versions
 of the LB and BP tests in the case of ARMA$(0,1)(0,1)_s$ model \eqref{ARMA21MonteCarlo}-\eqref{bruitARCH}  with $\alpha_1=0.45$.
The number of replications is $N=1000$. }}
\begin{center}
\begin{tabular}{c c ccc ccc c}
\hline\hline \\
$s$& Length $n$ & Lag $m$ & $\mathrm{{LB}}_{\textsc{sn}}$&$\mathrm{BP}_{\textsc{sn}}$&$\mathrm{{LB}}_{\textsc{w}}$&$\mathrm{BP}_{\textsc{w}}$&$\mathrm{{LB}}_{\textsc{s}}$&$\mathrm{BP}_{\textsc{s}}$
\vspace*{0.2cm}\\\hline
&& $4$&89.2 & 89.1 & 98.5 & 98.5& 100.0& 100.0\\
&& $8$&83.0 & 82.9 & 96.8 & 96.7 &100.0& 100.0\\
4 &$n=500$& $12$&71.5 & 71.5 & 96.3 & 96.3& 100.0& 100.0\\
 &&$15$&64.9&  64.7 & 95.6 & 95.6& 100.0& 100.0\\
 &&$18$&56.2 & 55.7 & 95.1 & 95.0& 100.0& 100.0\\
 &&$20$&50.6 & 49.4 & 94.5&  94.4& 100.0& 100.0\\
  \cline{2-9}
&& $4$&99.5 & 99.5& 100.0& 100.0& 100.0& 100.0\\
&& $8$&99.4 & 99.4& 100.0& 100.0& 100.0& 100.0\\
4 &$n=2,000$& $12$&99.1 & 99.1& 100.0& 100.0& 100.0& 100.0\\
 &&$15$&99.2 & 99.2 & 99.9 & 99.9& 100.0& 100.0\\
 &&$18$&99.5&  99.5& 100.0 &100.0& 100.0& 100.0\\
 &&$20$&99.2 & 99.2 & 99.9 & 99.9& 100.0 &100.0\\

\hline
&& $4$&90.1 & 90.1 & 98.4 & 98.3& 100.0& 100.0\\
&& $8$&82.4&  82.2&  97.3 & 97.2& 100.0& 100.0\\
12 &$n=500$& $12$&74.3 & 74.0&  97.0&  96.9& 100.0 &100.0\\
 &&$15$&68.7 & 68.3 & 96.9 & 96.5& 100.0 &100.0\\
 &&$18$&56.9 & 56.0 & 96.1&  95.7& 100.0& 100.0\\
 &&$20$&48.9 & 47.9 & 96.0  &95.7 &100.0 &100.0\\
 \cline{2-9}
&& $4$&99.6 & 99.6& 100.0& 100.0& 100.0& 100.0\\
&& $8$&99.6 & 99.6 & 99.9 & 99.9& 100.0& 100.0\\
12 &$n=2,000$& $12$&99.5 & 99.5& 100.0& 100.0& 100.0& 100.0\\
 &&$15$&99.4 & 99.4& 100.0& 100.0& 100.0& 100.0\\
 &&$18$&99.5 & 99.5& 100.0& 100.0& 100.0& 100.0\\
 &&$20$&99.0 & 99.0& 100.0&  99.9& 100.0 &100.0\\
\hline\hline
\\
\end{tabular}
\end{center}
\label{tabGARCHpuissance}
\end{table}
\subsection{Application to real data}
We now consider an application to monthly mean total sunspot number obtained by
taking a simple arithmetic mean of the daily total sunspot number over all days of each calendar month.
The observations (sunspot) covered the period from January 01, 2010 to December 31, 2018 which
correspond to $n=108$ observations. The series exhibit seasonal behavior $(s=12)$.
The data were obtain from the website of the World Data Center, Solar Influences Data Analysis Center,
Royal Observatory of Belgium (http://www.sidc.be/silso/datafiles).

\begin{figure}[hbt]
\vspace*{12cm} \protect \includegraphics{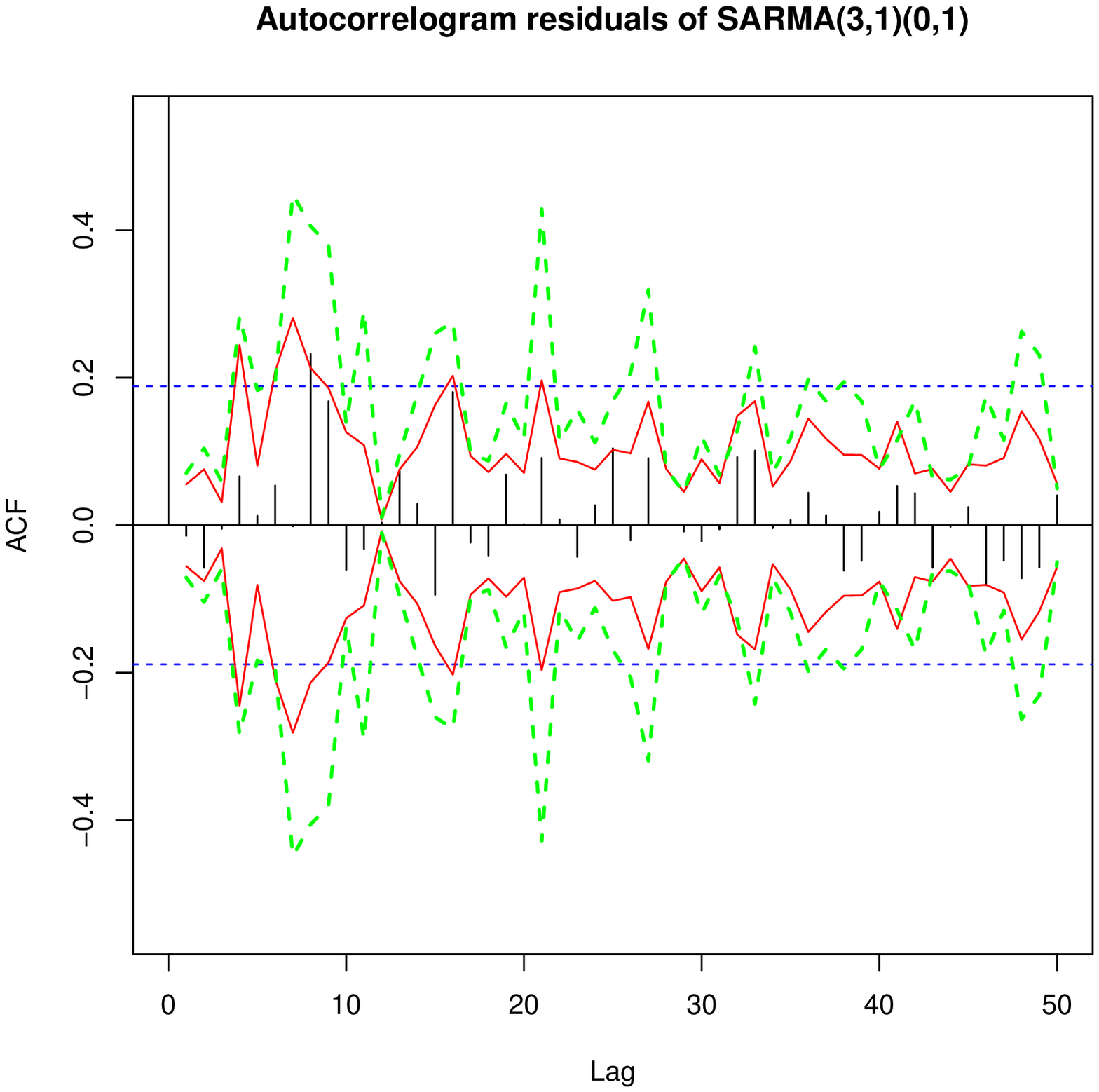} \vspace*{2.5cm}
\caption{\label{fig2} {\footnotesize Autocorrelation  of the particular SARMA$(3,1)(0,1)_{12}$ model
residuals for the  mean-logarithm-corrected of monthly mean total sunspot number.
The horizontal dotted lines (blue color) correspond to the 5\% significant limits obtained under the strong SARMA assumption.
The solid lines (red color) and  dashed lines (green color) correspond also  to the 5\% significant limits under the weak SARMA assumption.
The full lines  correspond to the asymptotic significance limits for the residual autocorrelations
obtained in Proposition~\ref{propo1}. The dashed lines (green color) correspond to the self-normalized asymptotic significance limits for the residual autocorrelations
obtained in Theorem~\ref{sn2}.}}
\end{figure}

Let $Z_t=\log(\mbox{sunspot}_t)-\log(\mbox{sunspot}_{t-1})$ and denoting by $X_t=Z_t-\mathbb{E}(Z_t)$ the mean-corrected
series. We adjust the particular SARMA$(3,1)(0,1)_{12}$ model of the form
$$X_{t}=a_{03}X_{t-3}+\epsilon_{t}-b_{01}\epsilon_{t-1}-\mathbf{b}_{01}\epsilon_{t-s}+b_{01}\mathbf{b}_{01}\epsilon_{t-s-1}.$$
The quasi-maximum likelihood estimators of $\theta_0=(a_{03},\mathbf{b}_{01},\mathbf{b}_{01})'$ were obtained as
\begin{eqnarray*}
\hat{\theta}_n=\left(\begin{array}{ccc}-0.1810 &[0.0253]&(0.0000)\\ 0.5438 &[0.0204]&(0.0000)\\-0.1139&[0.0237]& (0.0000)\end{array}\right)\text{ and }
\hat\sigma_{\epsilon}^2=0.1728,
\end{eqnarray*}
where the estimated asymptotic standard errors obtained from \eqref{NAsymp} (respectively the
$p$-values), of the estimated parameters (first column), are given into brackets
(respectively in parentheses).
We apply portmanteau tests to the residuals of this model.
Figure~\ref{fig2} displays the residual autocorrelations and their 5\% significance limits under the strong SARMA$(3,1)(0,1)_{12}$ and weak SARMA$(3,1)(0,1)_{12}$ assumptions.
In view of Figure~\ref{fig2}, the diagnostic
checking of residuals does not indicate any inadequacy.
All of the sample autocorrelations should lie between the bands (at 95\%)
shown as dashed lines (green color), solid lines (red color) and the horizontal dotted (blue color).
\section{Conclusion}\label{conclusion}
From these simulation experiments and from the asymptotic theory, we
draw the conclusion that the standard methodology, based on the
QMLE, allows to fit SARMA representations of a wide class of
nonlinear time series. But it is often restrictive to consider
that the innovation process is directly observed.
In future works, we intent to study how the existing estimation (see \cite{FZ2004,FZ2010}) and diagnostic checking (see \cite{Zhu2013}) procedures should be
adapted in  the situation where the GARCH process used in  these simulation experiments is not directly observed, but constitutes the
innovation of an observed SARMA-(seasonal)GARCH process which will be able to extend considerably the range of applications.
\clearpage
\appendix
\section{Proofs}\label{proof}

The proofs of Theorems \ref{sn1} and \ref{sn2} follow the same lines  as in  \cite{BMS2018} and are similar.
To have its own autonomy, the proofs will be rewrite and adapt.

\noindent{\bf Proof of Theorem \ref{sn1}}

We recall that the  Skorokhod space $\mathbb D^k [0{,}1]$ is the set of $\mathbb R^k-$valued functions defined on $[0{,}1]$ which are right continuous and have left limits.
It is endowed with the Skorokhod topology and the weak convergence on $\mathbb D^k [0{,}1]$ is mentioned by $\xrightarrow[]{\mathbb{D}^k}$. We finally denote by $\lfloor x \rfloor$
the integer part of the real $x$.

To prove the result we need to recall some following results of \cite{BM2014,frz,fz98,fz00,McLeod1978}.

We denote by $a_{i}^*$, $b_{i}^*$, $\mathbf{a}_{i}^*$ and $\mathbf{b}_{i}^*$ the coefficients defined by
\begin{eqnarray*}
a^{-1}_{\theta_0}(z)=\sum_{i\geq 0}a_{i}^*z^i,\;b^{-1}_{\theta_0}(z)=\sum_{i\geq 0}b_{i}^*z^i,\;\mathbf{a}^{-1}_{\theta_0}(z)=\sum_{i\geq 0}\mathbf{a}_{i}^*z^i \text{ and }\mathbf{b}^{-1}_{\theta_0}(z)=\sum_{i\geq 0}\mathbf{b}_{i}^*z^i,\quad |z|\leq 1.
\end{eqnarray*}
Following  \cite{frz} and \cite{McLeod1978} (see also \cite{BM2014}), the noise derivatives involving in the expression of $J(\theta_0)$ and $I(\theta_0)$ can be represented as
\begin{eqnarray}\label{serieLamda}
\frac{\partial \epsilon_t}{\partial\theta}&=&\sum_{i\geq 1}\lambda_i\epsilon_{t-i},\\ \lambda_i&=&\left(-a_{i-1}^*,\dots,-a_{i-p}^*,b_{i-1}^*,\dots,b_{i-q}^*,-\mathbf{a}_{i-1s}^*,
\dots,-\mathbf{a}_{i-Ps}^*,\mathbf{b}_{i-1s}^*,\dots,\mathbf{b}_{i-Qs}^*\right)'\in\mathbb{R}^{k_0},\nonumber
\end{eqnarray}
with $a_{i}^*=b_{i}^*=\mathbf{a}_{i}^*=\mathbf{b}_{i}^*=0$ when $i<0$.

For any $\theta\in\Theta\subset \mathbb R^{k_0}$ and any $(l,m)\in\{ 1,\dots,k_0\}^2$, under the above Assumptions, there exists absolutely summable and deterministic sequences  $(c_i(\theta))_{i\ge 0}$,  $(\lambda_{i,l}(\theta))_{i\ge 1}$  and $(\lambda_{i,l,m}(\theta))_{i\ge 1}$ such that, almost surely,
\begin{align}
\label{ce1}
\epsilon_t(\theta)&=\sum_{i=0}^{\infty}c_i(\theta)\epsilon_{t-i}\ ,\  \frac{\partial \epsilon_t(\theta)}{\partial \theta_l} =\sum_{i=1}^{\infty}\lambda_{i,l}(\theta)\epsilon_{t-i} \  \text{and}\   \frac{\partial^2\epsilon_t(\theta)}{\partial \theta_l\partial\theta_m} =\sum_{i=2}^{\infty}\lambda_{i,l,m}(\theta)\epsilon_{t-i}
\\
e_t(\theta)&=\sum_{i=0}^{t-1}c_i(\theta)e_{t-i}\ ,\  \frac{\partial e_t(\theta)}{\partial \theta_l} =\sum_{i=1}^{t-1}\lambda_{i,l}(\theta)e_{t-i} \  \text{and}\   \frac{\partial^2 e_t(\theta)}{\partial \theta_l\partial\theta_m} =\sum_{i=2}^{t-1}\lambda_{i,l,m}(\theta)e_{t-i}
\label{ce2}
\end{align}
with $c_0(\theta)=1$. A useful property of the above three sequences that they are asymptotically exponentially small. Indeed there exists $\rho\in ]0{,}1[$ and a positive constant $K$ such that, for all $i\ge 1$, we have
\begin{align}\label{rho}
\sup_{\theta\in \Theta} \Big ( |c_i(\theta)| +|\lambda_{i,l}(\theta)| + |\lambda_{i,l,m}(\theta)|\Big ) & \le K\, \rho^i \ .
\end{align}
See Lemmas A.1. and A.2. of \cite{fz00}  for a more detailed treatment.

Now, in view of \eqref{tayl} it is clear that the asymptotic behaviour of $\hat\gamma_m$ is related to the limit distribution of $w_t=\left (- \Upsilon_t'J^{-1'} {,} \epsilon_t\epsilon_{t-1} {,}\dots{,} \epsilon_t\epsilon_{t-m} \right )'$.
First, we prove that
$\frac{1}{\sqrt{n}}\sum_{j=1}^{\lfloor nr\rfloor }\Lambda w_j$ converges on the Skorokhod space to a Brownian motion.
More precesily, we have to show that
\begin{equation}
\label{fclt}
\frac{1}{\sqrt{n}}\sum_{j=1}^{\lfloor nr\rfloor }\Lambda w_j
\xrightarrow[n\to\infty]{\mathbb D^m}
\left(\Psi\Psi'\right)^{1/2} B_{m}(r)
\end{equation}
where $(B_m(r))_{r\ge 0}$ is a $m$-dimensional standard Brownian motion.

Using \eqref{serieLamda}, the process $w_t$ can be rewritten as
$$ w_t =  \left ( -2\sigma_0^{-2}\left\{\sum_{i=1}^{\infty} \lambda_{i,1}(\theta_0)\epsilon_t \epsilon_{t-i}\ , \ \cdots \ ,\ \sum_{i=1}^{\infty} \lambda_{i,k_0}(\theta_0)\epsilon_t \epsilon_{t-i}\right\}'J^{-1'}\  {,}\ \epsilon_t\epsilon_{t-1} \ {,}\ \dots{,} \ \epsilon_t\epsilon_{t-m}  \right )'$$
and thus the non-correlation between $\epsilon_t$ implies that $w_t$ has zero expectation with values in $\mathbb R^{k_0+m}$. In order to apply the functional central limit theorem for strongly mixing process,
we need to identify the asymptotic covariance matrix in the classical central limit theorem for the sequence $(w_t)_{t\ge 1}$.
It is proved in Subsection \ref{diagnostic1} that
\begin{align}\label{convvv}
\frac{1}{\sqrt{n}}\sum_{t=1}^{n}w_t \xrightarrow[n\to\infty]{\mathrm{d}}\mathcal{N}(0,\Xi:=2\pi f_w(0))
\end{align}
where $f_w(0)$ is the spectral density of the
stationary process $(w_t)_{t\in\mathbb Z}$ evaluated at frequency 0.  The main issue is to prove the existence of the matrix $\Xi$
which is a consequence {\bf (A3)} and Davydov's inequality  \cite{davy}.
For that sake, one has to introduce for any integer $k$, the random variables
\begin{align*}
w_t^k =\left ( -2\sigma_0^{-2}\left\{\sum_{i=1}^{k} \lambda_{i,1}(\theta_0)\epsilon_t \epsilon_{t-i}\ , \ \cdots \ ,\ \sum_{i=1}^{k} \lambda_{i,k_0}(\theta_0)\epsilon_t \epsilon_{t-i}\right\}'J^{-1'}\  {,}\ \epsilon_t\epsilon_{t-1} \ {,}\ \dots{,} \ \epsilon_t\epsilon_{t-m}  \right )'.
\end{align*}
Since $w^k$ depends on a finite number of values of the noise-process $\epsilon$, it also satisfies a mixing property (see Theorem 14.1 in \cite{davidson1994}, p. 210).
Based on the Davydov inequality  (see \cite{davy}), the arguments developed in the Lemma A.1 in \cite{frz} (see also \cite{fz98}) imply that
\begin{align}\label{convvvk}
\frac{1}{\sqrt{n}}\sum_{t=1}^{n}w^k_t \xrightarrow[n\to\infty]{\mathrm{d}}\mathcal{N}(0,\Xi_k)
\end{align}
where
\begin{equation*}
\Xi_k:=2\pi f_{w^k}
(0)=\sum_{h=-\infty}^{+\infty}\cov(w^k_t,w^k_{t-h})=
\sum_{h=-\infty}^{+\infty}\e (w_t^k {{w}^{k}_{t-h}}' )
\end{equation*}
and thus \eqref{convvv} holds. Moreover we have that $\lim_{k\to\infty}\Xi_k =\Xi$.

Since the matrix $\Xi$ is positive definite, it can be factorized as $\Xi=\Delta\Delta'$ where the $(k_0+m)\times (k_0+m)$ lower triangular matrix $\Delta$ has nonnegative diagonal entries. Therefore, we have
$$\frac{1}{\sqrt{n}}\sum_{t=1}^{n}\Lambda w_t \xrightarrow[n\to\infty]{\mathrm{d}}\mathcal{N}(0,\Lambda\Xi\Lambda'),$$
and the new variance matrix can also been factorized as $\Lambda\Xi\Lambda'= (\Lambda\Delta)(\Lambda\Delta)':=\Psi\Psi'$, where $\Psi\in\mathbb{R}^{m\times (k_0+m)}$.  Thus,
${n}^{-1/2}\sum_{t=1}^{n}\left(\Psi\Psi'\right)^{-1/2}\Lambda w_t \xrightarrow[n\to\infty]{\mathrm{d}} \mathcal{N}\left(0,I_{m}\right)$ where $I_{m}$ is the identity matrix of order $m$.
The above arguments also apply to matrix $\Xi_k$ with some matrix $\Psi_k$ which is defined analogously as $\Psi$. Consequently,
$$\frac{1}{\sqrt{n}}\sum_{t=1}^{n}\Lambda w_t^k \xrightarrow[n\to\infty]{\mathrm{d}}\mathcal{N}(0,\Lambda\Xi_k\Lambda'), $$
and we also have ${n}^{-1/2}\sum_{t=1}^{n}\left(\Psi_k\Psi_k'\right)^{-1/2}\Lambda w^k_t \xrightarrow[n\to\infty]{\mathrm{d}} \mathcal{N}\left(0,I_{m}\right)$

Now we are able to apply the functional central limit theorem for strongly mixing process of \cite{herr}. We have for any $r\in (0{,}1)$,
\begin{align*}
\frac{1}{\sqrt{n}}\sum_{j=1}^{\lfloor nr\rfloor }\left(\Psi_k\Psi_k'\right)^{-1/2}\Lambda w_j^k
\xrightarrow[n\to\infty]{\mathbb D^m}
B_{m}(r) .
\end{align*}
For all $j\in\{1,\dots,[nr]\}$, we write $$\left(\Psi\Psi'\right)^{-1/2}\Lambda w_j^k = \Big (\left(\Psi\Psi'\right)^{-1/2}-\left(\Psi_k\Psi_k'\right)^{-1/2}\Big)\Lambda w_j^k  +  \left(\Psi_k\Psi_k'\right)^{-1/2}\Lambda w_j^k$$
and we obtain that
\begin{align*}
\frac{1}{\sqrt{n}}\sum_{j=1}^{\lfloor nr\rfloor }\left(\Psi\Psi'\right)^{-1/2}\Lambda w_j^k
\xrightarrow[n\to\infty]{\mathbb D^m}
B_{m}(r) .
\end{align*}
In order to conclude \eqref{fclt}, it remains to observe that, uniformly with respect to $n$,
\begin{align}
\label{ll}
Z^k_n(r):=\frac{1}{\sqrt{n}}\sum_{j=1}^{\lfloor nr\rfloor }\left(\Psi\Psi'\right)^{-1/2}\Lambda Y_j^k
\xrightarrow[k\to\infty]{\mathbb D^m}
0 ,
\end{align}
where
$$ Y_t^k = \left ( -2\sigma_0^{-2}\left\{\sum_{i=k+1}^{\infty} \lambda_{i,1}(\theta_0)\epsilon_t \epsilon_{t-i}\ , \ \cdots \ ,\ \sum_{i=k+1}^{\infty}
\lambda_{i,k_0}(\theta_0)\epsilon_t \epsilon_{t-i}\right\}'J^{-1'}\  {,}\ \epsilon_t\epsilon_{t-1} \ {,}\ \dots{,} \ \epsilon_t\epsilon_{t-m}
 \right )' \ .$$
By Lemma 4 in \cite{fz98}, we have
$$
\sup_{n} \mathrm{Var} \left ( \frac{1}{\sqrt{n}} \sum_{j=1}^n Y_j^k \right ) \xrightarrow[k\to\infty]{} 0
$$
and since $\lfloor nr\rfloor\le n$,
$$\sup_{0\le r\le 1}\sup_n \left (  \| Z_n^k(r) \| \right ) \xrightarrow[k\to\infty]{} 0 . $$
Thus \eqref{ll} is true and the proof of the first step \eqref{fclt} is achieved.

The previous step ensures us that  Assumption 1 in \cite{lobato} is satisfied for the sequence $(\Lambda w_t)_{t\ge 1}$.
We follow the arguments developed in Sections 2 and 3 in \cite{lobato}, the second step is to show that
\begin{align}
\label{conv}
C_m=\frac{1}{n^2}\sum_{t=1}^{n}S_{t}S'_{t}   \xrightarrow[n\to\infty]{\mathrm{d}} \left(\Psi\Psi'\right)^{1/2}{V}_{m}\left(\Psi\Psi'\right)^{1/2},
\end{align}
by applying the continuous mapping theorem on the Skorokhod space and
where the random variable $V_m$ is defined in \eqref{UkVk}. The main issue is to obtain that
\begin{align}\nonumber
\frac{1}{\sqrt{n}}{S}_{[nr]} &=\frac{1}{\sqrt{n}}\sum_{j=1}^{[nr]}\left({\Lambda}{w}_j-\Lambda\bar{w}\right)
=\frac{1}{\sqrt{n}}\sum_{j=1}^{[nr]}{\Lambda}{w}_j-\frac{[nr]}{n}\left(\frac{1}{\sqrt{n}}\sum_{j=1}^{n}\Lambda {w}_j\right)\\
&\xrightarrow[n\to\infty]{\mathbb D^m}
\left(\Psi\Psi'\right)^{1/2} B_{m}(r)-r\left(\Psi\Psi'\right)^{1/2} B_{m}(1),\label{convS}
\end{align}
by continuous mapping theorem and using \eqref{fclt}, the fact that $[nr]/n\to r$ as $n\to\infty$.
In view of \eqref{convS}, it follows that
\begin{align*}
C_m&=\frac{1}{n^2}\sum_{t=1}^{n}S_{t}S'_{t} =\frac{1}{n}\sum_{t=1}^{n}\int_{t/n}^{(t+1)/n}S_{[nr]}S'_{[nr]}dr
=\int_{1/n}^{(n+1)/n}\left(\frac{1}{\sqrt{n}}{S}_{[nr]}\right)\left(\frac{1}{\sqrt{n}}{S}_{[nr]}\right)'dr
\\&  \xrightarrow[n\to\infty]{\mathrm{d}} \left(\Psi\Psi'\right)^{1/2}\left(\int_0^1\left ( {B}_{m}(r)-r{B}_{m}(1)\right)\left({B}_{m}(r)-r{B}_{m}(1)\right)'dr\right)\left(\Psi\Psi'\right)^{1/2}= \left(\Psi\Psi'\right)^{1/2}{V}_{m}\left(\Psi\Psi'\right)^{1/2},
\end{align*}
which prove \eqref{conv}.
Since $\sqrt{n} \hat\gamma_m = n^{-1/2}\sum_{t=1}^n\Lambda w_t+\mathrm{o}_\mathbb{P}(1)$, using (\ref{fclt}) and (\ref{conv})  we obtain
\begin{align*}
n\hat\gamma_m'C_m^{-1}\hat\gamma_m & = \left(\frac{1}{\sqrt{n}}\sum_{t=1}^{n}\Lambda w_t+\mathrm{o}_\mathbb{P}(1) \right)'C_{m}^{-1}\left(\frac{1}{\sqrt{n}}\sum_{t=1}^{n}\Lambda w_t+\mathrm{o}_\mathbb{P}(1) \right) \\
 &\xrightarrow[n\to\infty]{\mathrm{d}} \left(\left[\Psi\Psi'\right]^{1/2}{B}_{m}(1)\right)'\left(\left[\Psi\Psi'\right]^{1/2}{V}_{m}\left\{\left[\Psi\Psi'\right]^{1/2}\right\}'\right)^{-1}\left(\left[\Psi\Psi'\right]^{1/2}{B}_{m}(1)\right)
\\&\qquad={B}'_{m}(1){V}_{m}^{-1}{B}_{m}(1)=:\mathcal U_m.
\end{align*}
The proof of Theorem \ref{sn1} is then complete.\zak
\noindent{\bf Proof of Theorem \ref{sn2}:}

We write $\hat C_m=C_m+ \nabla_n$ where $ \nabla_n  ={n}^{-2} \sum_{t=1}^n \big (S_tS_t' - \hat S_t\hat S_t' \big)$.
There are three kinds of entries in the matrix $\nabla_n$. The first one is a sum composed of
$$\upsilon_t^{k,k'} = \epsilon_t^2(\theta_0)\epsilon_{t-k}(\theta_0)\epsilon_{t-k'}(\theta_0)  - {e}_t^2(\hat\theta_n){e}_{t-k}(\hat\theta_n){e}_{t-k'}(\hat\theta_n)  $$
for $(k,k')\in\{1,\dots,m\}^2$.
Using \eqref{rho} and the consistency of $\hat\theta_n$, we have $\upsilon_t^{k,k'} = \mathrm{o}(1)$ almost surely.
The two last  kinds of entries of $\nabla_n$ come from the following quantities for $i,j\in\{1,\dots,k_0\}$ and $k\in\{1,\dots,m\}$
\begin{align*}
\tilde\upsilon_t^{k,i} & = \epsilon_t^2(\theta_0)\epsilon_{t-k}(\theta_0) \frac{\partial \epsilon_{t}(\theta_0)}{\partial \theta_i}-  {e}_t^2(\hat\theta_n){e}_{t-k}(\hat\theta_n)\frac{\partial {e}_{t}(\hat\theta_n)}{\partial \theta_i},
\\
\bar\upsilon^{i,j}_t & = \epsilon^2_{t}(\theta_0)\frac{\partial \epsilon_{t}(\theta_0)}{\partial \theta_i} \frac{\partial \epsilon_{t}(\theta_0)}{\partial \theta_j}-{e}^2_{t}(\hat\theta_n)\frac{\partial {e}_{t}(\hat\theta_n)}{\partial \theta_i}
\frac{\partial {e}_{t}(\hat\theta_n)}{\partial \theta_j}
\end{align*}
and they also satisfy $\tilde\upsilon_t^{k,i} + \bar{\upsilon}_t^{i,j}= \mathrm{o}(1)$ almost surely by using \eqref{serieLamda} and \eqref{rho}.
Consequently, $\nabla_n= \mathrm{o}(1)$ almost surely as $n$ goes to infinity.
Thus one may find a matrix $\nabla^\ast_n$, that tends to the null matrix almost surely, such  that
\begin{align*}
n\, \hat\gamma_m'\hat C_m^{-1}\hat\gamma_m &  = n\, \hat\gamma_m' (C_m + \nabla_n)^{-1}\hat\gamma_m   = n\, \hat\gamma_m' C_m^{-1}\hat\gamma_m + n\, \hat\gamma_m'\nabla_n^\ast \hat\gamma_m \ .
\end{align*}
Thanks to the arguments developed in the proof of Theorem \ref{sn1}, $n \hat\gamma_m' C_m^{-1}\hat\gamma_m $ converges in distribution. So $n \hat\gamma_m'\nabla_n^\ast \hat\gamma_m$ tends to zero in distribution, hence in probability.  Then $n \hat\gamma_m'\hat C_m^{-1}\hat\gamma_m $ and  $n \hat\gamma_m' C_m^{-1}\hat\gamma_m $ have the same limit in distribution and the result is proved. \zak
\begin{center}
{\bf Acknowledgements}

 We sincerely thank the anonymous reviewers and Editor for helpful remarks.
The authors wish to acknowledge the support from the "Séries temporelles et valeurs extrêmes : théorie et 
applications en modélisation et estimation des risques" Projet Région (Bourgogne Franche-Comté, France) grant No OPE-2017-0068.
\end{center}

\bibliographystyle{plain}
\bibliography{biblio-yac}


\end{document}